\documentclass{amsart}[12pt]
\usepackage{amsmath, amsthm, amscd, amsfonts}

\makeatletter \oddsidemargin.9375in \evensidemargin \oddsidemargin
\theoremstyle{definition}

\theoremstyle{remark}

\numberwithin{equation}{section}

\begin{document}

%\noindent {}\\[1.00in]

\title[Orthogonality of
Homogeneous geodesics  ...]{Orthogonality of Homogeneous
geodesics  on the tangent bundle }

\author[R.Chavosh Khatamy]{R.Chavosh Khatamy $^{*}$\\}

\subjclass[2000]{53C30, 53C22.}

\keywords{ Hadamard matrices,tangent bundle,
   homogeneous manifold, solvable Lie group,
  homogeneous (geodesic)
vector.\\
$^{*}$\textit{Department of Mathematics, Islamic Azad University
Tabriz Branch,
\\Tabriz-IRAN\\}
E-mail: \texttt{chavosh@tabrizu.ac.ir, r\_chavosh@iaut.ac.ir}
\indent }

\begin{abstract}

Let  $\xi=(G\times_{K} \mathcal{G} / \mathcal{K}, \rho_{\xi},
\emph{G} / \emph{K},\mathcal{G} / \mathcal{K})$  be the associated
bundle and $\tau_{G/K}=(T_{G/K},\pi_{G/K},G/K, \textrm{R}^{m})$ be
the tangent bundle of special examples of odd dimension solvable
Lie groups equipped with left invariant Riemannian metric. In this
paper we  prove some  conditions about the existence of
homogeneous geodesic  on the base space of $\tau_{G/K}$ and
homogeneous (geodesic) vectors on the fiber space of $\xi$ .

\end{abstract} \maketitle

\section{Introduction and preliminaries }

Let  $G$ be a connected Lie group and $K$ be a closed subgroup of
$G$. The set of left cosets of $K$ in $G$ is denoted by $G/K$ and
can be given a unique differentiable structure ([6], vol.II,
chap.2), and hence $M=G/K$ is called a homogeneous manifold. When
a Lie group $G$ acts transitively isometric on a Riemannian
manifold $M$, we can identify $M$ with the set $G/K$ of left
cosets of the isotropy group $K$ of a point $x_{0}\in M$. The
point $x_{0}$ is called the origin of $M$. Let $\bigtriangledown$
be an affine connection on $M=G/K$ and let $\bigtriangledown$ be
invariant under the natural action of $T:G\times M \longrightarrow
M$. Then a geodesic $\gamma :I \longrightarrow M$ is called a {\em
homogeneous geodesic} if, there exists a 1-parameter subgroup
$t\longrightarrow \exp tX$ , $t\in \textbf{R},$  of $G$ with $X\in
\mathcal{G}$$=T_eG$
 such that
$$\gamma(t)=T(\exp tX, x_{0}).$$
 Where $\gamma(0)=x_{0} \in M$,
and $exp:\mathcal{G} \rightarrow \emph{G} $ is the exponential map
[9].\\

\textbf{Definition 1.1.} A vector $0\neq X\in \mathcal{G}$ is
called a \textit{homogeneous vector} (or \emph{geodesic vector}),
if the curve $\gamma(t)=(\exp tX)(x_0)$ is a geodesic on $M=G/K$
[9].

The following result can be found in [7], proposition 1.

 \emph{Any homogeneous Riemannian manifold $G/K$ has the reductive
 decomposition of the form}

$$\mathcal{G}=\mathcal{M}+\mathcal{K}$$

\emph{where $\mathcal{M} \subset \mathcal{G}$ is a vector
subspace, such that $Ad(K)(\mathcal{M})\subset \mathcal{M}$.}

Let $M=G/K$  be a  Riemannian manifold and
$\mathcal{G}=\mathcal{M}+\mathcal{K}$, its reductive
decomposition. Then the natural map $\phi:G\longrightarrow G/K=M$
induces a linear epimorphism $(d\phi)_e:T_eG\longrightarrow
T_{x_0}M$, and the vector space $\mathcal{M}$  can be identified
with $T_{x_0}M$.
 If $C$ is a scalar product on $\mathcal{M}$ induced by the
scalar product on $T_{x_0}M$, then the following lemma holds (see
[9], proposition 2.1).\\

\textbf{Lemma 1.2.} {\em If $X$ belongs to $\mathcal{G}$, let
$[X,Y]_{\mathcal{M}}$ and $X_{\mathcal{M}}$  be the components of
$[X,Y]$ and $X$ in $\mathcal{M}$ with respect to reductive
decomposition, then $X$ is homogeneous vector} (or \emph{geodesic
vector}) \emph{iff}

$$C(X_\mathcal{M},[X,Y]_{\mathcal{M}})=0 \quad \forall Y \in
{\mathcal{G}}.$$\\

\textbf{Proposition 1.3.}([10]). \emph{A finite family $\{\gamma
_{1},\gamma_{2},\ldots ,\gamma_{n} \}$ of homogeneous geodesics
through $x_{o}\in M$ is orthogonal ( respectively, linearly
independent) if the $\mathcal{ M}$-component of the corresponding
homogeneous vectors are orthogonal (respectively, linearly
independent)}.

Let $\wp=(P,\pi,B,G)$ be a smooth fiber bundle. A pair $(\wp,T)$
is called a (smooth)\emph{principal bundle with structure group
$G$}, if $T:P\times G\longmapsto P$ is a right action of $G$ on
$P$ and $\wp$ admits a coordinate representation
$\{(U_{\alpha},\psi_{\alpha})\}$ such that \\
$$\psi_{\alpha}(x,ab)=\psi_{\alpha}(x,a)b,\hspace{0.5cm}  x\in
U_{\alpha},\hspace{0.5cm} a,b\in G, $$(see [6], vol.II, chap.V). \\
Let $\wp=(P,\pi,B,G)$ be a principal bundle and $F$ be a
differentiable manifold. Consider the left action $Q$, of $G$ on
the product manifold $P\times F$ given by
$$Q_{a}(z,y)=(z,y)a=(za,a^{-1}y)\hspace{0.5cm} z\in P,
y\in F, a\in G .$$\\
The set of orbits of this action is denoted by $P\times_{G} F$ and
$$ q: P\times F\rightarrow P\times_{G} F $$ will denote the
corresponding projection, i.e.,  $q(z,y)$ is the orbit through
$(z,y)$. The map  $q$ determines a map $\rho_{\xi}: P \times_{G} F
\rightarrow B$ such that,
$$\rho_{\xi} \circ q =\pi \circ \pi_{p} . $$
Where, $ \pi_{p}:P\times F\rightarrow P$ is the canonical
projection and $\pi:P\rightarrow B$ is the bundle map.\\
There is a unique smooth structure on $P\times_{G} F$, such that
$\xi=(P\times_{G} F,\rho_{\xi},B,F)$ is a smooth fiber bundle (see [6], vol.II, chap.V, sec.2).\\

\textbf{Definition 1.4.} The fiber bundle $\xi=(P\times_{G}
F,\rho_{\xi},B,F)$, is called the
 \emph{associated bundle} with  $\wp=(P,\pi,B,G).$\\
Let $K$ be a closed subgroup of $G$. The principal fiber bundle
$\Im=(G,\pi,G/K,K)$, is called
\emph{homogeneous bundle,} (See [3]).\\

Let $\Im=(G,\pi,G/K,K)$ be a fiber bundle with group structure
$K$, and  let $G$ be a connected Lie group and $K$ a closed
subgroup of $G$, (see [1], definition 2.2). We take the Lie
algebras $\mathcal{G}$ and $\mathcal{K}$ of $G$ and $K$
respectively, in [1] and [2], we proved some relations between the
homogeneous vector in the fiber space of the associated bundle,
$\xi=(G\times_{K} \mathcal{G} / \mathcal{K}, \rho_{\xi}, \emph{G}
/\emph{K},\mathcal{G} / \mathcal{ K})$  and the homogeneous
geodesic in the base space of a principal homogeneous bundle
$\Im=(G,\pi,G/K,K)$.In [3], we consider the homogeneous bundle
$\Im=(G,\pi,G/K,K)$ and the tangent bundle $\tau_{G/K}$ of
$M=G/K$, and  give some results about the existence of homogeneous
vectors on  the fiber space of $\tau_{G/K}$, for both cases of $G$
semisimple and weakly semisimple Lie group.

Now, we investigate the existence of mutually orthogonal linearly
independent homogeneous geodesics in the base space of the tangent
bundle  $\tau_{G}$ of homogeneous Riemannian manifold $G$ given in
theorem 2.2.

\section{Main results}

Let $\Im=(G,\pi,G/K,K)$ be a principal homogeneous bundle,
 with the associated bundle  $\xi=(G\times_{K}
\mathcal{ K} , \rho_{\xi}, \emph{G} / \emph{K}, \mathcal{ G} /
\mathcal{K})$. Let

$G$ be the matrix group of all matrices of the form

$$ \left ( \begin {array}{cccccc} e^{z_{0}} & 0 & \ldots & 0 & x_{0} \\ 0 & e^{z_{1}} & \ldots &
0 & x_{1} \\ \ldots & \ldots & \ldots & \ldots & \ldots \\ 0 & 0 &
\ldots & e^{z_{n}} & x_{n}
\\ 0 & 0 & \ldots  & 0 & 1 \end {array} \right )  $$
where, $(x_{0}, x_{1}, \cdots, x_{n}, z_{1} \cdots, z_{n})\in
\textbf{R}^{2n+1}$. The Lie group $G$ is unimodular and solvable
(see [8], pp.134-136), with the  left invariant Riemannian metric

$$g=\sum_{i=0}^{n} e^{-2z_{i}}dx_{i}^{2}+ \lambda^{2}\sum_{k,j=0}^{n} dz_{k}dz_{j}.$$

Where $\lambda \neq 0$ is a constant. Then $G$ is a homogeneous
Riemannian manifold with the origin at $(0, 0, \cdots, 0)$ ([8],
p.134).

 Let $\mathcal{G}=\mathcal{M}+\mathcal{K}$ be the reductive
decomposition of $\mathcal{G}$ , then $\mathcal{K}$= 0, and hence
$\mathcal{G}=\mathcal{M}$.

In [3], we prove the following lemma\\

\textbf{Lemma 2.1.}\emph{ Let $\Im=(G,\pi,G/K,K)$, be a
homogeneous bundle. Then
 $$\xi=(G\times_{K} \mathcal{G} / \mathcal{K},
\rho_{\xi}, \emph{G} / \emph{K},\mathcal{G} / \mathcal{K}),$$ is
the
associated bundle of $\Im=(G,\pi,G/K,K)$.}\\

By lemma 2.1, we can take $\xi=(G\times \mathcal{M} , \rho_{\xi},
\emph{G} , \mathcal{M})$,
be the associated bundle of $\Im=(G,\pi,G/K,K)$.\\

In [4], we let $G$ be a 3-dimensional solvable Lie group, given in
[8], pp.134, and prove some  results about the existence of
homogeneous  vectors on the fiber
 space of $\tau_{G/K}$ and $\xi$.\\
In [5], we extend theorems 5.6 and  5.7 in [4], and give the
following theorem, for the odd dimensional solvable Lie group.\\

\textbf{Theorem 2.2.}([5]). \emph{Let $\Im=(G,\pi,G/K,K)$, be a
principal homogeneous bundle and } $\xi=(G\times_{K} \mathcal{G} /
\mathcal{K}, \rho_{\xi}, \emph{G} / \emph{K}, \mathcal{G} /
\mathcal{K}),$ \emph{be the associated bundle of
$\Im=(G,\pi,G/K,K)$. If $G$ is the matrix group of all matrices of
the form}

$$ \left ( \begin {array}{cccccc} e^{z_{0}} & 0 & \ldots & 0 & x_{0} \\ 0 & e^{z_{1}} & \ldots &
0 & x_{1} \\ \ldots & \ldots & \ldots & \ldots & \ldots \\ 0 & 0 &

\ldots & e^{z_{n} } &  x_{n}
\\ 0 & 0 & \ldots  & 0 & 1 \end {array} \right )  $$

\emph{ where $(x_{0}, x_{1}, \cdots, x_{n}, z_{1} \cdots,
z_{n})\in \textbf{R}^{2n+1}$ and $z_{0}=-(z_{1}+z_{2}+\cdots
+z_{n})$. Then a vector $V$ in the fiber space of $\xi$ is a
homogenous (geodesic ) if and only if its components
$$(x_{0}, x_{1}, \cdots, x_{n}, z_{1}, \cdots, z_{n})$$
 satisfy the following conditions}

$$ x_{0}(z_{1}+z_{2}+\cdots +z_{n})=0
\hspace{1cm} x_{1}z_{1}=0, \cdots, x_{n}z_{1}=0$$
$$x_{0}^{2}-x_{1}^{2}=0, \cdots, x_{0}^{2}-x_{n}^{2}=0.$$

 In the proof of the Theorem 5.3 in [3], we give a strong
isomorphism between  the tangent bundle
$$\tau_{G/K}=(T_{G/K},\pi_{G/K},G/K, \textbf{R}^{m})$$
 and the associated bundle
 $$\xi=(G\times_{K} \mathcal{ G} / \mathcal{K},\rho_{\xi}, \emph{G} / \emph{K},\mathcal{G} / \mathcal{K}),$$
 then  under hypothesis of theorem 2.2  there is a strong isomorphism between,the associated bundle
$$\xi=(G\times \mathcal{M} , \rho_{\xi}, \emph{G} , \mathcal{M})$$
and the tangent bundle
$$\tau_{G}=(T_{G},\pi_{G},G, \textbf{R}^{2n+1})$$
so we have,\\

\textbf{Corollary 2.3.}([5]). \emph{With hypothesis of theorem
2.2, let}
$$\tau_{G}=(T_{G},\pi_{G},G, \textbf{R}^{2n+1})$$
\emph{be the tangent bundle of the homogeneous Riemannian manifold
$G$. Then a vector $W$ in $ \textbf{R}^{2n+1}$ is a homogeneous
vector (under isomorphism), if and only if its component $(x_{0},
x_{1}, \cdots, x_{n}, z_{1} \cdots, z_{n})$  satisfy the following
conditions}

$$ x_{0}(z_{1}+z_{2}+\cdots +z_{n})=0
\hspace{1cm} x_{1}z_{1}=0, \cdots, x_{n}z_{1}=0$$
$$x_{0}^{2}-x_{1}^{2}=0, \cdots, x_{0}^{2}-x_{n}^{2}=0.$$

 In [3], theorem 5.4. we  give a subspace of $\mathcal{G'}$ such that all
member of this subspace are homogeneous vectors, and by strong
isomorphism between $\tau_{G/K}$ and $\xi$ we can find a subspace
of $\textbf{R}^m$ (under isomorphism) such that all members of
this subspace are homogeneous vectors,

In the following theorem, we consider the tangent bundle

  $$\tau_{G}=(T_{G},\pi_{G},G, \textbf{R}^{2n+1})$$

of the homogeneous Riemannian manifold $G$ in theorem 2.2, and
give structure of all subspaces  of $\textbf{R}^{2n+1}$ such that
all their  members are homogeneous vectors.\\

\textbf{Theorem 2.4.}([5]). \emph{ Let}
$$\tau_{G}=(T_{G},\pi_{G},G, \textbf{R}^{2n+1})$$
\emph{be the tangent bundle of the homogeneous Riemannian manifold
$G$, (given in theorem 2.2). Then all  homogeneous vectors are
decomposed into an n-dimension vector subspace $W$ in
$\textbf{R}^{2n+1}$ and $2^{n}$, one-dimension vector subspace in
$\textbf{R}^{2n+1}$ generated by all vectors of the form }
$X_{0}\pm X_{1}\pm \cdots \pm X_{n}.$\\

 By proposition 1.3 and Theorem 2.4 we have, the following
result  about linearly independence of homogeneous geodesics on
 the base space of $\tau_{G}$\\

\textbf{Corollary 2.5.}([5]). \emph{With hypothesis of theorem
2.2, the tangent bundle  }
$$\tau_{G}=(T_{G},\pi_{G},G, \textbf{R}^{2n+1})$$
\emph{admits $2n+1$ linearly independent homogeneous geodesics
through the origin $\{e\}$ of the base space of $\tau_{G}$.}\\

Now, we investigate orthogonality of homogeneous vectors on the
fiber space of tangent bundle,
$$\tau_{G}=(T_{G},\pi_{G},G, \textbf{R}^{2n+1}).$$
In [3] we prove some conditions about existence and orthogonality
of homogeneous vectors for both cases of $G$ semisimple and weakly
semisimple. For example in theorem 5.3 in [3], we prove that if
$G$ is a semisimple Lie group then there are $m$  orthogonal
homogeneous vectors on the fiber space of the tangent bundle,
 $$\tau_{G/K}=(T_{G/K},\pi_{G/K},G/K,
\textrm{R}^{m})$$

In the follow, we want to get some conditions about linearly
independent and orthogonality of homogeneous vectors on the fiber
space and homogeneous geodesics on the base space of the tangent
bundle of the homogeneous Riemannian manifold $G$ (given in
theorem 2.2). For this we need to considering to relations between
orthogonality of homogeneous vectors and the  Hadamard matrices.\\

\textbf{Definition 2.6.} A \emph{Hadamard matrix of order $k$ }is
$k\times k$ square matrix whose entries are all equal to $\pm1$,
and such that $A.A^{t}=kI_{k}$, where $I_{k}$ is the unit
matrix.\\

The condition   $A.A^{t}=kI_{k}$, in definition 2.6, implies that
the $k$ rows or columns of a Hadamard matrix represent orthogonal
$k$-tuples, with all entries equal to +1 or -1, we can use this
fact for considering to structure of Hadamard matrices and
orthogonality of homogeneous (geodesic) vectors. \\

\textbf{Lemma 2.7.} \emph{Let $\tau_{G}=(T_{G},\pi_{G},G,
\textbf{R}^{2n+1})$ be the tangent bundle of the homogeneous
Riemannian Lie group $G$ of all matrices of the form}

$$ \left ( \begin {array}{cccccc} e^{z_{0}} & 0 & \ldots & 0 & x_{0} \\ 0 & e^{z_{1}} & \ldots &
0 & x_{1} \\ \ldots & \ldots & \ldots & \ldots & \ldots \\ 0 & 0 &

\ldots & e^{z_{n}} &  x_{n}
\\ 0 & 0 & \ldots  & 0 & 1 \end {array} \right )  $$

\emph{ where $(x_{0}, x_{1}, \cdots, x_{n}, z_{1} \cdots,
z_{n})\in \textbf{R}^{2n+1}$ and $z_{0}=-(z_{1}+z_{2}+\cdots
+z_{n})$ ,}\\
then;\\
 (i)\emph{ If $(n+1)$ is odd, then there are not any two mutually
 orthogonal $(n+1)$-tuples with all entries equal to $\pm 1$.}\\
(ii)\emph{If $(n+1)$ is even and not  divisible by $4$, then there
are exactly two mutually
orthogonal $(n+1)$-tuples with all entries equal to $\pm 1$.}\\

\textbf{Proof.} Let $\tau_{G}=(T_{G},\pi_{G},G,
\textbf{R}^{2n+1})$ be the tangent bundle of the homogeneous
Riemannian manifold $G$ (given in theorem 2.2). By  Corollary 2.3,
and theorem 2.4 a vector $w$ in $ \textbf{R}^{2n+1}$ is a
homogeneous (geodesics) vector (under isomorphism), if and only if\\
$$ A)\hspace{1cm} w \in W= span (Z_{1}, Z_{2}, \cdots, Z_{n})$$
$$ B)\hspace{1cm} w=\sum_{i=0}^{i=n}x_{i}X_{i}\hspace{0.5cm} and \hspace{0.5cm}x_{0}^{2}-x_{1}^{2}=0, \cdots,
x_{0}^{2}-x_{n}^{2}=0.$$\\

 As concerns homogeneous
(geodesics) vectors of type (B), they are all generated by the
vectors of the form $X_{0}+\epsilon_{1} X_{1}+ \cdots
+\epsilon_{n} X_{n},$ where $\epsilon_{i}\in \{1,-1\}$. Therefore,
the problem of finding mutually orthogonal geodesics vectors of
type (B) is equivalent to the algebraic problem of finding
$(n+1)$-tuples, with all entries equal to $\pm1$, which are
mutually orthogonal with respect to the standard
scaler product in $\textbf{R}^{n+1}$.\\
Let $(n+1)$ be odd number and $W_{1}$ and $W_{2}$ be two
$(n+1)$-tuples with all entries equal to $\pm 1$. The scaler
product of $W_{1}$ and $W_{2}$ is the sum of the products of their
entries and all such products are equal to $\pm 1$. By hypotheses,
$(n+1)$ is odd, then sum of the products of their entries dose not
vanish, so $W_{1}$ and $W_{2}$ can not be orthogonal, so we obtain
$(i)$. For the second statement of the lemma, we spouse that
$(n+1)= 2m$, where $m$ is odd, let $V_{1}$ and $V_{2}$ be two
$(n+1)$-tuples with all entries equal to $\pm 1$, such that
$V_{1}= ( 1, 1, \cdots, 1)$ and $V_{2}= (-1, 1, \cdots -1, 1)$,
then $V_{1}$ and $V_{2}$ are orthogonal. Now, we spouse that $V$,
$W$, $Z$, are three mutually orthogonal $(n+1)$-tuples with all
entries equal to $\pm 1$. Then, we compute the  scaler product of
$V$, $W$ and  $Z$ by $V$. In this way, we can obtain three
mutually orthogonal $(n+1)$-tuples vectors $V'$, $W'$, $Z'$ such
that all entries equal to $\pm 1$. If we take $V'=(-1, -1, \cdots
,-1)$, then by orthogonality of $V'$ and $W'$, $W'$ has exactly
$m$ entries equal to $-1$ and exactly $m$ entries equal to $1$. We
then multiply, component by component, and applying a fixed
permutation of the all  entries for mutually orthogonal
$(n+1)$-tuples vectors $V'$, $W'$, $Z'$, such that this
applications will preserve the orthogonality of $V'$, $W'$, $Z'$.
By this way, we can obtain $W'= (1, 1, \cdots, 1, -1, -1, \cdots,
-1)$, but $m$ is odd and the orthogonality of $V'$, $W'$, $Z'$ is
imposable, this gives a contradiction, and the proof of the lemma
is complete. $\Box$\\

Before starting some additional results, we recall the fact that
$A.A^{t}=kI_{k}$, in definition 2.6 implies that the $k$ rows or
columns of a Hadamard matrix represent orthogonal $k$-tuples, with
all entries equal to +1 or -1, for the case $n+1$ be divisible by
4, the problem related to algebraic problem of the existence of
Hadamard matrices of order $n+1$. Therefore, we get at once the
following proposition.\\

 \textbf{Proposition 2.8.}{\emph{With hypothesis of lemma
2.7, let $n+1$ be divisible by 4, then $\textbf{R}^{2n+1}$ admits
$n+1$ mutually orthogonal $(n+1)$-tuples vectors with all entries
equal to $\pm 1$, if and only if, there exists a Hadamard matrices
of
order $n+1$.}\\

Now, we can prove the following theorem about the linearly
independent and the maximum number of the orthogonal homogeneous
(geodesic) vectors on the fiber space of

$$\tau_{G}=(T_{G},\pi_{G},G, \textbf{R}^{2n+1}).$$\\

\textbf{Theorem 2.9.} \emph{Let $\tau_{G}=(T_{G},\pi_{G},G,
\textbf{R}^{2n+1})$ be the tangent bundle of the homogeneous
Riemannian Lie group $G$, (given in theorem 2.2 and lemma 2.7)
then;}\\

(i)\emph{There are $2n+1$ linearly independent homogeneous
(geodesics) vectors in the fiber space of through the
$\tau_{G}$.}\\
(ii)\emph{ The maximum number of the orthogonal homogeneous
(geodesic) vectors on the fiber space of $\tau_{G}$ is $n+1$, in
the case that $n+1$ is odd. }\\
(ii)\emph{The maximum number of the orthogonal homogeneous
(geodesic) vectors on the fiber space of $\tau_{G}$, is $n+2$, in
the case that $n+1$ is even and  not divisible by $4$. }\\
(iv)\emph{The maximum number of the orthogonal homogeneous
(geodesic) vectors on the fiber space of $\tau_{G}$, is $2n+1$, in
the case that $n+1$ is even and divisible by $4$ and there exists
a Hadamard matrices of order $n+1$.}\\

\textbf{Proof.} Theorem 2.4 and corollary 2.3, conclude the fist
part of theorem, it is easy to see that, there exist $n+1$
linearly independent homogeneous (geodesics) vectors of type (B),
( see proof of lemma 2.7), then there are $2n+1$ linearly
independent homogeneous (geodesics) vectors in the fiber space of $\tau_{G}$.}\\
The second and the third part of the theorem follows from (i) and
(ii), in lemma 2.7. Finally, as an immediate consequence from
proposition 2.8, we
obtain (iv).$\square$\\

By proposition 1.3 and theorem 2.9 we complete corollary 2.5 about
the number of
  linearly independent homogeneous geodesics through origin of the
  base space of $\tau_{G}$.\\

\textbf{Corollary 2.10.}  \emph{Let $\tau_{G}=(T_{G},\pi_{G},G,
\textbf{R}^{2n+1})$ be the tangent bundle of the homogeneous
Riemannian Lie group $G$, (given in theorem 2.2 and lemma 2.7)
then;}\\

(i)\emph{There are $2n+1$ linearly independent homogeneous
geodesics vectors through the origin $\{e\}$ of the base space of $\tau_{G}$.}\\
(ii)\emph{ The maximum number of the orthogonal homogeneous
geodesic through the origin $\{e\}$ of the base space of
$\tau_{G}$, is $n+1$, in
the case that $n+1$ is odd. }\\
(ii)\emph{The maximum number of the orthogonal homogeneous
geodesic through the origin $\{e\}$ of the base space of
$\tau_{G}$, is $n+2$, in
the case that $n+1$ is even and  not divisible by $4$. }\\
(iv)\emph{The maximum number of the orthogonal homogeneous
geodesic  through the origin $\{e\}$ of the base space of
$\tau_{G}$, is $2n+1$, in the case that $n+1$ is even and
divisible by $4$ and there exists a Hadamard matrices of order
$n+1$.}\\

{\bf Acknowledgement}\\
The author would like to express his appreciation of professor O.
Kowalski for his invaluable suggestions and also professor M.
Toomanian for his constructive comments.

The author was supported by the funds of the Islamic Azad
University- Tabriz Branch, (IAUT).

\end{document}